\documentclass{article}
\usepackage{amsmath,amssymb,amsthm}
\usepackage{graphicx,subfigure,float,url,color}
\usepackage{mathrsfs}
\usepackage[colorlinks=true]{hyperref}
\usepackage{pdfsync}
\usepackage{enumitem}

\def\R{\mathbb{R}}

\renewcommand{\geq}{\geqslant}
\renewcommand{\leq}{\leqslant}

\newtheorem{theorem}{Theorem}  
\newtheorem{proposition}{Proposition}

\theoremstyle{definition}\newtheorem{remark}{Remark}

\title{Linear quadratic optimal control turnpike in finite and infinite dimension: two-term expansion of the value function}

\author{
Veljko A\v{s}kovi\'c\footnote{Sorbonne Universit\'e, CNRS, Universit\'e Paris Cit\'e, Inria, Laboratoire Jacques-Louis Lions (LJLL), F-75005 Paris, France (\texttt{veljkoaskovic@hotmail.com}).}
\and
Emmanuel Tr\'elat\footnote{Sorbonne Universit\'e, CNRS, Universit\'e Paris Cit\'e, Inria, Laboratoire Jacques-Louis Lions (LJLL), F-75005 Paris, France (\texttt{emmanuel.trelat@sorbonne-universite.fr}).}
\and
Hasnaa Zidani\footnote{Insa Rouen Normandie, Normandie Université, Laboratoire de Mathématiques (LMI),  UR 3226, F-76000 Rouen, France (\texttt{hasnaa.zidani@insa-rouen.fr}).}
}

\date{}

\begin{document}

\maketitle

\begin{abstract}
In this paper, we consider a linear quadratic (LQ) optimal control problem in both finite and infinite dimensions. We derive an asymptotic expansion of the value function as the fixed time horizon $T$ tends to infinity. The leading term in this expansion, proportional to $T$, corresponds to the optimal value attained through the classical turnpike theory in the associated static problem. The remaining terms are associated with optimal stabilization problems towards the turnpike.
\end{abstract}




\section{Introduction and Formulation of the problem}\label{sec1}

\subsection{Setting}

The optimization of linear autonomous control systems has been a subject of considerable interest in the field of control theory.  In this work, we are concerned with a long-time optimal control problem. We aim to analyze the connection between this problem and a static problem with an infinite horizon. Our study concerns both the finite-dimensional and infinite-dimensional cases.

Consider two   Hilbert spaces, denoted as  $(X,\langle\ ,\ \rangle_X)$ and $(U,\langle\ ,\ \rangle_U)$, each equipped with their respective scalar products.  Our focus lies on the guidance of a linear control system from an initial point to a final point within a fixed time interval $T>0$. Specifically, we fix the values of $y_0, y_1, y_d \in X$, and $u_d \in U$, the linear quadratic (LQ) optimal control problem consists of steering the  system
\begin{equation}\label{contsyst}
\dot y(t) = A y(t) + B u(t) 
\end{equation}
from the initial point $y(0)=y_0$ to the final point $y(T)=y_1$, with a control input  $u\in L^2([0,T],U)$ that minimizes the following cost functional
\begin{equation}\label{cost}
C_T(u) = \frac{1}{2} \int_0^T \left( \Vert y(t)-y_d\Vert_Q^2 + \Vert u(t)-u_d\Vert_R^2 \right) dt .
\end{equation}
Here, we use the notations
$\Vert y-y_d\Vert_Q^2 = \langle Q(y-y_d),y-y_d\rangle_X$
and
$\Vert u-u_d\Vert_R^2 = \langle R(u-u_d),u-u_d\rangle_U$.

For a finite-dimensional case, where $X = \R^n$ and $U = \R^m$ (for some  integers  $n\geq 1$ and $m\geq 1$),  the system is defined by a $n\times n$ matrix $A$,  and a control matrix $B$ of dimension $n\times m$.  The matrices $Q$ and $R$ that appear in the cost function are symmetric positive definite   of sizes 
$n\times n$ and $m\times m$, respectively.  All the matrices are real-valued. 
 
In the general case (notably, in infinite dimension),   $A:D(A)\rightarrow X$ is a linear operator on the Hilbert space $X$ of domain $D(A)$ generating on $X$ a $C_0$ semi-group $(e^{tA})_{t\geq 0}$ and 
$B\in L(U,X)$ is a linear bounded operator from $U$ to $X$ (see \cite{CurtainZwart,trelat_book,TucsnakWeiss,Zabczyk} for the general framework).
The operators  $Q\in L(X)$ and $R\in L(U)$ are positive definite, boundedly invertible, selfadjoint operators respectively on $X$ and $U$. 

To ensure the well-posedness of the optimal control problem, we introduce the following assumption \ref{H2} to guarantee exact controllability in the state space X.
\begin{enumerate}[label=$\bf(H)$]
\item\label{H2} 
There exists $T_0>0$ such that the control system \eqref{contsyst} is exactly controllable in the state space $X$, with controls $u\in L^2([0,T],U)$, in any time $T>T_0$.
\end{enumerate}
Note that, in finite-dimensional spaces,
\ref{H2} is equivalent to the Kalman condition on the pair $(A,B)$ and does not depend on $T$ (and one can take $T_0=0$). 
In infinite-dimensional spaces, 
the exact controllability condition \ref{H2} is satisfied for instance for wave equations in appropriate functional spaces under some geometric conditions (with a minimal controllability time required) but it does not hold for heat equations, although such equations enjoy approximate controllability properties.


Now, under Assumption \ref{H2} and by strict convexity of the cost function,  the optimal control problem \eqref{contsyst}-\eqref{cost} has a unique solution for every $T>T_0$ (see, e.g., \cite{LeeMarkus,trelat_book,Troltzsch,Zabczyk}), denoted $(y_T(\cdot),u_T(\cdot))$.
The central focus of this paper is on the \emph{value function} of the above optimal control problem, defined by
\begin{equation}\label{def_valuefunction}
V_T(y_0,y_1) = \min \{ C_T(u)\ \mid\ y(0)=y_0, \ y(T)=y_1 \} = C_T(u_T) ,
\end{equation}
that is the minimal cost required to steer the control system \eqref{contsyst} from $y_0$ to $y_1$ in time $T$. 
The objective is to provide a two-term expansion of $V_T(y_0,y_1)$ for large time $T$ and to identify the various components of this expansion. 

In the subsequent sections, we explore the turnpike property and the main term of the asymptotic expansion. The turnpike property, well-established in the literature, suggests that optimal trajectories exhibit a specific structure, involving rapid transitions to and from a turnpike state. 

\subsection{Turnpike and main term of the asymptotic expansion}\label{sec_turnpike}
Determining an asymptotic value to $V_T(y_0, y_1)$ as $T\rightarrow+\infty$
 proves to be a straightforward task,  closely related to the renowned \emph{turnpike property}, as succinctly outlined in \cite{TrelatZhangZuazua_SICON2018, TrelatZuazua_JDE2015}. This task involves using the underlying connection between the value function $V$ and the dynamic behavior of the system over prolonged time horizons, aligning with the well-established principles represented by the turnpike property.

\subsubsection{Static optimal control problem}\label{sec_static}
Let $(\bar y,\bar u)\in D(A)\times U$ 
be the unique solution of the (strictly convex) constrained optimization problem
given by
\begin{equation}\label{static_problem}
\min_{Ay+Bu=0} \frac{1}{2}\left( \Vert y-y_d\Vert_Q^2 + \Vert u-u_d\Vert_R^2 \right) 
\end{equation}
referred to as the \emph{static optimal control problem}.
According to the Karush-Kuhn-Tucker (KKT, see \cite{Kurcyusz,Troltzsch}) rule, there exists $\bar\lambda\in D(A^*)$ such that
\begin{equation}\label{KKT}
Q(\bar y-y_d) = A^*\bar\lambda 
\qquad\textrm{and}\qquad
R(\bar u-u_d) = B^*\bar\lambda .
\end{equation}
Here, we have used the adjoint operators $A^*:D(A^*)\rightarrow X$ and $B^*:U\rightarrow X$,  and we have identified $X'\simeq X$ and $U'\simeq U$.
%
As a consequence of \eqref{KKT}, we have $\bar y=y_d+Q^{-1}A^*\bar\lambda$ and $\bar u=u_d+R^{-1}B^*\bar\lambda$ and thus
\begin{equation}\label{syst_KKT}
\underbrace{
\begin{pmatrix}
A & BR^{-1}B^* \\
Q & -A^*
\end{pmatrix}
}_M
\begin{pmatrix}\bar y\\ \bar\lambda\end{pmatrix}
=
\begin{pmatrix}-Bu_d\\ Qy_d\end{pmatrix},
\end{equation} 
where $M$ is a linear operator on $X\times X$ of domain $D(M)=D(A)\times D(A^*)$.
We denote by
\begin{equation}\label{Vbar}
\bar V = \frac{1}{2}\left( \Vert \bar y-y_d\Vert_Q^2 + \Vert \bar u-u_d\Vert_R^2 \right) 
\end{equation}
the optimal value of the optimization problem \eqref{static_problem}. 
Note that
$$
\bar V = -\frac{1}{2} \langle Ay_d+Bu_d,\bar\lambda\rangle_X . 
$$

\subsubsection{Application of the Pontryagin maximum principle}\label{sec_PMP}
For every $T>T_0$, by the Pontryagin maximum principle (PMP, see \cite{LeeMarkus,LiYong,Pontryagin,trelat_book,Troltzsch}) applied to the optimal control problem \eqref{contsyst}-\eqref{cost}, of optimal solution $(y_T(\cdot),u_T(\cdot))$, there exists an (unique) absolutely continuous costate $\lambda_T(\cdot):[0,T]\rightarrow D(A^*)$ satisfying almost everywhere on $[0,T]$
\begin{equation}\label{eq_adj}
\dot\lambda_T(t) = -A^*\lambda_T(t) + Q(y_T(t)-y_d) ,
\end{equation} 
and we have $u_T(t) = u_d + R^{-1}B^*\lambda_T(t)$ for almost every $t\in[0,T]$. 

The non-triviality of the Pontryagin maximum principle is not always guaranteed when $\dim X=+\infty$ and may require stringent conditions \cite{Casas97,RayZid98}.  Indeed, in the infinite-dimensional case, it is well-known that the PMP may fail (see \cite{LiYong}) if there is an infinite number of constraints on the terminal states, which is the case here since the initial and final states are prescribed.

Actually,  Assumption \ref{H2} (exact controllability) implies that the controllability Gramian operator is an isomorphism (see \cite{Flandoli}), and then the Hilbert Uniqueness Method (HUM), see \cite{Lions_SIREV}), can be applied. However, in the LQ case, HUM exactly coincides with the PMP (see \cite{trelat_book}). Another justification is that, under \ref{H2}, the differential Riccati theory can be applied (see \cite{Flandoli,Zabczyk}), which leads as well to the adjoint equation \eqref{eq_adj}. In the proof of the main result, we will indeed revisit these issues and particularly focus on the Riccati theory. 

\subsubsection{Exponential turnpike property.  Main term of the asymptotic expansion}

It has been proven in \cite{TrelatZuazua_JDE2015} for finite dimensional problems  and in \cite{TrelatZhangZuazua_SICON2018} for infinite dimension that, under the assumption \ref{H2}, there exist $C,\nu>0$ such that
\begin{equation}\label{exp_turnpike}
\Vert y_T(t)-\bar y\Vert_X + \Vert u_T(t)-\bar u\Vert_U + \Vert \lambda_T(t)-\bar\lambda\Vert_X \leq C e^{-\nu t(T-t)} \quad \forall t\in[0,T] \quad\forall T>T_0.
\end{equation}
The inequality \eqref{exp_turnpike} is referred to as the \emph{exponential turnpike property}. It implies that, except at the beginning and the end of the time frame $[0,T]$, the "dynamic" optimal triple $(y_T(\cdot),u_T(\cdot),\lambda_T(\cdot))$ is exponentially close to the "static" optimal triple $(\bar y,\bar u,\bar\lambda)$.
It's worth noting that the constants $C$ and $\nu$ are independent of $T>T_0$.
 This property ensures that the optimal solution remains close to its static counterpart over the majority of the time interval, providing a stable and predictable behavior for the optimal trajectory.
%
As an immediate consequence of this exponential turnpike property, it follows that
$$
V_T(y_0,y_1) = T \bar V + \mathrm{o}(T)
$$
as $T\rightarrow +\infty$, where $V_T(y_0,y_1)$ is defined by \eqref{def_valuefunction} and $\bar V$ is defined by \eqref{Vbar}.

In the following sections, we present our main result, which provides the second term in the large-time asymptotic expansion of $V_T(y_0, y_1)$, that is, an equivalent of $V_T(y_0, y_1) - T \bar V$ as $T\rightarrow +\infty$.

\section{Two-term asymptotic expansion of the value function.  Main result}\label{sec_main_result}
\subsection{A preliminary comment}

The turnpike property was identified in the 1950s by Nobel Prize laureate Samuelson and his coauthors in \cite{DorfmanSamuelsonSolow}, primarily in the context of econometrics (for historical insights and a comprehensive bibliography, refer to \cite{TrelatZuazua_JDE2015}). In essence, the turnpike property stipulates that the optimal trajectory approximately consists of three segments: the first (resp. the third) arc is short and represents a rapid transition to (resp., from) the turnpike, while the second, middle arc is long and consists of remaining at an optimal steady-state that is the turnpike. Notably, in the 1970s and 1980s, an equivalent property known as the exponential dichotomy property was identified in \cite{AndersonKokotovic,WildeKokotovic}. In that work, the authors demonstrated that large-time optimal trajectories can be approximated by concatenating two infinite-time trajectories. Each of these trajectories is the solution to an optimal stabilization (towards the turnpike) problem and serves as an approximation of the initial or final transient arc.

This is why, unsurprisingly, we proceed to preliminarily define two optimal stabilization problems. The first one aims to stabilize the forward-in-time control system \eqref{contsyst} from the (initial) point $y_0$ to the turnpike point $\bar y$, by minimizing a quadratic cost that measures the discrepancy between the trajectory and the turnpike. The second one aims to stabilize the backward-in-time control system \eqref{contsyst} from the (final) point $y_1$ to the turnpike point $\bar y$, again by minimizing the same quadratic cost.   To achieve this, we rely on the following assumption, which we consider satisfied throughout this section.

\begin{enumerate}[label=$\bf(H_A)$]
\item\label{HA}
The operator $A$ is the infinitesimal generator on $X$ of a $C_0$-group $(e^{tA})_{t\in \R}$.
\end{enumerate}

A well-known necessary and sufficient condition for \ref{HA} to hold is that both $A$ and $-A$ generate a $C_0$ semi-group.
It is worth noting that, like \ref{H2}, this assumption is satisfied for wave equations but not for heat equations.

\subsection{Forward stabilization problem} \label{sec_forward_stab}
Considering the forward control system with initial condition
$$
\dot y(t) = Ay(t)+Bu(t), \qquad y(0)=y_0,
$$
we define
$$
S_f(y_0) = \inf_{u\in L^2([0,+\infty),U)} \frac{1}{2} \int_0^{+\infty} \left( \Vert y(t)-\bar y\Vert_Q^2 + \Vert u(t)-\bar u\Vert_R^2 \right) dt .
$$
Under \ref{H2}, this infinite time horizon optimal control problem, called the \emph{forward stabilization problem}, has a unique solution $(y_f(\cdot),u_f(\cdot))$. Equivalently, setting $z=y-\bar y$ and $v=u-\bar u$, we have
$$
\dot z(t) = Az(t)+Bv(t), \qquad z(0)=y_0-\bar y,
$$
and
$$
S_f(y_0) = \inf_{v\in L^2([0,+\infty),U)} \frac{1}{2} \int_0^{+\infty} \left( \Vert z(t)\Vert_Q^2 + \Vert v(t)\Vert_R^2 \right) dt .
$$
Since, by \ref{H2}, the pair $(A,B)$ is stabilizable, by the well known Riccati algebraic theory (see \cite{CurtainZwart, LeeMarkus, Zabczyk}), the latter ``shifted" stabilization problem has a unique solution $(z_f(\cdot),v_f(\cdot))$ with the feedback control
$v_f = -R^{-1}B^*Pz_f$
where the linear bounded operator $P:X\rightarrow X$ is the unique nonnegative selfadjoint solution of the algebraic Riccati equation (of unknown $\mathtt{X}$)
\begin{equation}\label{ARE}
A^*\mathtt{X}+\mathtt{X}A-\mathtt{X}BR^{-1}B^*\mathtt{X} = -Q .
\end{equation}
Actually, $P$ is positive definite and is even boundedly invertible, due to the exact controllability assumption \ref{H2} and the group assumption \ref{HA} (see \cite{Flandoli}). 
Moreover, the operator
$$
A_- = A-BR^{-1}B^*P ,
$$
of domain $D(A_-)=D(A)$, generates on $X$ an exponentially stable $C_0$ group $(e^{tA_-})_{t\in R}$, and
$S_f(y_0) = \frac{1}{2} \langle Pz(0),z(0)\rangle_X$.
Therefore $u_f = \bar u - R^{-1}B^*P(y_f-\bar y)$ and $y_f(t)=\bar y + e^{tA_-}(y_0-\bar y)$, and
$$
S_f(y_0) = \frac{1}{2} \langle P(y_0-\bar y),y_0-\bar y\rangle_X .
$$

\subsection{Backward stabilization problem}\label{sec_backward_stab}
We now turn our attention to the backward control system  with initial condition
$$
\dot y(t) = -Ay(t)-Bu(t), \qquad y(0)=y_1,
$$
we define
$$
S_b(y_1) = \inf_{u\in L^2([0,+\infty),U)} \frac{1}{2} \int_0^{+\infty} \left( \Vert y(t)-\bar y\Vert_Q^2 + \Vert u(t)-\bar u\Vert_R^2 \right) dt .
$$
This infinite time horizon optimal control problem, called the \emph{backward stabilization problem}, has a unique solution $(y_b(\cdot),u_b(\cdot))$.  Here again,  setting $z=y-\bar y$ and $v=u-\bar u$, we obtain
$$
\dot z(t) = -Az(t)-Bv(t), \qquad z(0)=y_1-\bar y,
$$
and
$$
S_b(y_1) = \inf_{v\in L^2([0,+\infty),U)} \frac{1}{2} \int_0^{+\infty} \left( \Vert z(t)\Vert_Q^2 + \Vert v(t)\Vert_R^2 \right) dt .
$$
Under \ref{H2}, since $(e^{tA})_{t\in\R}$ is a group (\ref{HA} is crucially used here), the pair $(-A,-B)$ is stabilizable.

Here, let us emphasize an important observation. 
Denoting temporarily by $(R_f)$ the algebraic Riccati equation \eqref{ARE} associated with the pair $(A,B)$, the algebraic Riccati equation $(R_b)$ associated with the pair $(-A,-B)$ is $-A^*\mathtt{X}-\mathtt{X}A-\mathtt{X}BR^{-1}B^*\mathtt{X} = -Q$ (of unknown $\mathtt{X}$), and obviously $\mathtt{X}$ is a solution of $(R_f)$ if and only if $-\mathtt{X}$ is a solution of $(R_b)$.
According to  the Riccati algebraic theory, $(R_b)$ has a unique non-negative self-adjoint solution, denoted by $-N$, which is actually positive definite and even boundedly invertible. Therefore the linear boundedly invertible operator $N:X\rightarrow X$ is the unique non-positive self-adjoint solution of \eqref{ARE}, and actually is negative definite. 

Moreover, as per the Riccati algebraic theory, the ``shifted" backward stabilization problem possesses a unique solution $(z_b(\cdot),v_b(\cdot))$ with the feedback control
$v_b = -R^{-1}B^*N z_b$.
Defining the operator

$$
A_+ = A-BR^{-1}B^*N
$$
of domain $D(A_+)=D(A)$, the operator $-A_+$ (also of domain $D(A)$) generates an exponentially stable $C_0$ group $(e^{-tA_+})_{t\in\R}$, and
$S_b(y_1) = -\frac{1}{2} \langle Nz(0),z(0)\rangle_X$.
Therefore, $u_b = \bar u - R^{-1}B^*N(y_b-\bar y)$ and $y_b(t)=\bar y + e^{-tA_+}(y_1-\bar y)$, which yield
$$
S_b(y_1) = -\frac{1}{2} \langle N (y_1-\bar y),y_1-\bar y\rangle_X .
$$

\subsection{Main result}
Let $\nu$ be the exponential decay rate of $e^{tA_-}$ and of $e^{-tA_+}$ (they are actually the same).

\begin{theorem}\label{main_thm}
We have
\begin{equation}\label{asympt_exp}
V_T(y_0,y_1) = T\bar V + S_f(y_0) - \langle\bar\lambda,y_0-\bar y\rangle_X + S_b(y_1) + \langle\bar\lambda,y_1-\bar y\rangle_X + \mathrm{O}(e^{-\nu T})
\end{equation}
as $T\rightarrow+\infty$.
\end{theorem}

Additionally to the statement of Theorem \ref{main_thm}, we have the following results, stated in the next two propositions.

\begin{proposition}\label{prop_dissip}
Defining the optimal cost
\begin{equation}\label{def_Vf}
V_f(y_0) = \inf_{u\in L^2([0,+\infty),U)} \frac{1}{2} \int_0^{+\infty} \left( \Vert y(t)-y_d\Vert_Q^2 + \Vert u(t)-u_d\Vert_R^2 - \Vert \bar y-y_d\Vert_Q^2 - \Vert \bar u-u_d\Vert_R^2 \right) dt 
\end{equation}
for the forward control system,
and the optimal cost
\begin{equation}\label{def_Vb}
V_b(y_1) = \inf_{u\in L^2([0,+\infty),U)} \frac{1}{2} \int_0^{+\infty} \left( \Vert y(t)-y_d\Vert_Q^2 + \Vert u(t)-u_d\Vert_R^2 - \Vert \bar y-y_d\Vert_Q^2 - \Vert \bar u-u_d\Vert_R^2 \right) dt
\end{equation}
for the backward control system, we have
$$
V_f(y_0) = S_f(y_0-\bar y) - \langle\bar\lambda,y_0-\bar y\rangle_X,
\qquad
V_b(y_1) = S_f(y_0-\bar y) - \langle\bar\lambda,y_0-\bar y\rangle_X ,
$$
and \eqref{asympt_exp} is equivalent to
\begin{equation}\label{asympt_exp2}
V_T(y_0,y_1) = T\bar V + V_f(y_0) + V_b(y_1) + \mathrm{O}(e^{-\nu T})
\end{equation}
as $T\rightarrow+\infty$.
\end{proposition}

This alternative expression of the two-term large-time asymptotic expansion of $V_T(y_0,y_1)$ is interesting because the function inside the integral defining the cost in \eqref{def_Vf} and in \eqref{def_Vb} is exactly the supply rate function used to characterize the dissipativity property of the optimal control problem. The form \eqref{asympt_exp2} can thus be seen as a preliminary to a generalization to nonlinear optimal control problems, treated in \cite{ATZ}.

\begin{proposition}\label{prop_optsol}
In addition, as regards the optimal solution, we have
\begin{eqnarray}
y_T(t) &=& \bar y + z_f(t) + z_b(T-t) - e^{-(T-t)A_+} z_f(T)  - e^{tA_-} z_b(T) \nonumber\\ 
&& \phantom{\bar y + z_f(t) + z_b(T-t)} + e^{-(T-t)A_+} \mathrm{O}(e^{-2\nu T}) + e^{tA_-}  \mathrm{O}(e^{-2\nu T}) \label{exp_yT_1} \\
&=& \bar y + e^{tA_-}\left(y_0-\bar y-e^{-TA_+}(y_1-\bar y)+\mathrm{O}(e^{-2\nu T})\right) \nonumber\\
&& \phantom{\bar y} + e^{-(T-t)A_+}\left(y_1-\bar y+e^{TA_-}(y_0-\bar y)+\mathrm{O}(e^{-2\nu T}) \right) \label{exp_yT_2}
\end{eqnarray}
where
\begin{equation}\label{zfzb}
z_f(t) = e^{tA_-}(y_0-\bar y), \qquad z_b(t) = e^{-tA_+} (y_1-\bar y) ,
\end{equation}
and
\begin{eqnarray}
\lambda_T(t) &=& \bar\lambda - P z_f(t) - N z_b(T-t) + P e^{tA_-} z_b(T) - N e^{-(T-t)A_+} z_f(T) \nonumber\\
&& \phantom{\bar\lambda - P z_f(t) - N z_b(T-t)} + e^{-(T-t)A_+} \mathrm{O}(e^{-2\nu T}) + e^{tA_-}  \mathrm{O}(e^{-2\nu T}) \label{exp_lambdaT_1}\\
&=& \bar\lambda -P e^{tA_-}\left(y_0-\bar y-e^{-TA_+}(y_1-\bar y)+\mathrm{O}(e^{-2\nu T})\right) \nonumber\\
&& \phantom{\bar\lambda} - N e^{-(T-t)A_+}\left(y_1-\bar y-e^{TA_-}(y_0-\bar y)+\mathrm{O}(e^{-2\nu T})\right) \label{exp_lambdaT_2}
\end{eqnarray}
and, using that $\bar u = u_d+R^{-1}B^*\bar\lambda$,
\begin{eqnarray}
u_T(t) &=& u_d + R^{-1}B^*\lambda_T(t) \nonumber\\
&=& \bar u + v_f(t) + v_b(T-t) \label{exp_uT_1}\\
&=& \bar u - R^{-1}B^*P e^{tA_-}(y_0-\bar y) - R^{-1}B^*N e^{-(T-t)A_+}(y_1-\bar y) + \mathrm{O}(e^{-\nu T}) \label{exp_uT_2}
\end{eqnarray}
as $T\rightarrow+\infty$. 
\end{proposition}

The formulas \eqref{exp_yT_1}, \eqref{exp_lambdaT_1} and \eqref{exp_uT_1} above are interesting because they show how the optimal solution of \eqref{contsyst}-\eqref{cost} is related to the optimal solutions of the (shifted) forward and backward stabilization problems defined in Section \ref{sec_main_result}: $(z_f(\cdot),v_f(\cdot))$ (resp., $(z_b(\cdot),v_b(\cdot))$) is the optimal solution of the forward (resp., backward) stabilization problem defined in Section \ref{sec_forward_stab} (resp., Section \ref{sec_backward_stab}). Note that the corresponding costate, in the infinite time horizon PMP, is $\lambda_f(\cdot)=-Pz_f(\cdot)$ (resp., $\lambda_b(\cdot)=Nz_b(\cdot)$).

The formulas \eqref{exp_yT_2}, \eqref{exp_lambdaT_2} and \eqref{exp_uT_2} give the first terms of an expansion of the optimal solution within the scale given by $e^{TA_-}$ and $e^{-TA_+}$. The expansion can actually be obtained at any order (see the proof in Section \ref{sec_proof_exp}).
As a consequence, an expansion of the value function $V_T(y_0,y_1)$ can also be obtained at any order. In \eqref{asympt_exp} we have given only the two first terms of that expansion, because we know how to give an interpretation of those two terms; but it seems that terms of higher order do not have a nice interpretation.

\begin{remark}
It is a classical result of optimal control that the initial or final costates are given in terms of the gradient of the value function (sensitivity analysis), namely, $\lambda_T(0) = -\frac{\partial V_T}{\partial y_0}(y_0,y_1)$ and $\lambda_T(T) = \frac{\partial V_T}{\partial y_1}(y_0,y_1)$. It is then interesting to note that, thanks to \eqref{asympt_exp2} and \eqref{exp_lambdaT_2}, we have
$$
-\lambda_T(0) = \frac{\partial V_T}{\partial y_0}(y_0,y_1) = \frac{\partial V_f}{\partial y_0}(y_0,y_1)+\mathrm{O}(e^{-\nu T}) = -\bar\lambda+\underbrace{P(y_0-\bar y)}_{\frac{\partial S_f}{\partial y_0}(y_0)}+\mathrm{O}(e^{-\nu T}),
$$
and
$$
\lambda_T(T) = \frac{\partial V_T}{\partial y_1}(y_0,y_1) = \frac{\partial V_b}{\partial y_1}(y_1)+\mathrm{O}(e^{-\nu T}) = \bar\lambda\underbrace{-N(y_1-\bar y)}_{\frac{\partial S_b}{\partial y_1}(y_1)}+\mathrm{O}(e^{-\nu T}),
$$
as $T\rightarrow+\infty$. 
This is of interest, in particular, in view of initializing a numerical shooting method (see \cite{TrelatZuazua_JDE2015}).
\end{remark}

\section{Proof of Theorem \ref{main_thm} and of the subsequent remarks}
The proof goes in several steps, performed in the subsequent sections.
\subsection{Preliminary}
In Section \ref{sec_PMP}, we have applied the Pontryagin maximum principle, leading to the extremal system
\begin{equation}\label{extr_syst}
\begin{split}
\dot y_T &= Ay_T + BR^{-1}B^*\lambda_T+Bu_d \\
\dot\lambda_T &= Qy_T-A^*\lambda_T-Qy_d
\end{split}
\end{equation}
with $y(0)=y_0$ and $y(T)=y_1$. 
Noting that the pair $(\bar y,\bar\lambda)$ defined in Section \ref{sec_static} satisfies \eqref{syst_KKT}, setting
$$
\delta y_T(t) = y_T(t)-\bar y,\qquad \delta\lambda_T(t) = \lambda_T(t)-\bar\lambda\qquad \forall t\in[0,T] ,
$$
we get from \eqref{extr_syst} that
\begin{equation}\label{extr_syst_delta}
\frac{d}{dt}\begin{pmatrix} \delta y \\ \delta\lambda\end{pmatrix} = M \begin{pmatrix} \delta y \\ \delta\lambda\end{pmatrix}
\end{equation}
where we recall that $M=\begin{pmatrix}
A & BR^{-1}B^* \\
Q & -A^*
\end{pmatrix}$,
with $\delta y_T(0)=y_0-\bar y$ and $\delta y_T(T)=y_1-\bar y$. 

\subsection{Diagonalization by blocks of $M$}\label{sec_diago}
Let us first prove that $M$ is diagonalizable by blocks and is boundedly invertible. 
The following argument is borrowed from \cite{WildeKokotovic}; its generalization to infinite dimension is straightforward under the assumptions of exact controllability and of group.

In Section \ref{sec_main_result}, we have defined the boundedly invertible selfadjoint operators $P>0$ and $N<0$. 
Defining on $X\times X$ the linear bounded operator
$$
T = \begin{pmatrix}
\mathrm{id} & \mathrm{id} \\
-N & -P
\end{pmatrix} ,
$$
we first note that $T$ is boundedly invertible and
$$
T^{-1} = \begin{pmatrix}
\Delta^{-1} P & \Delta^{-1} \\
-\Delta^{-1} N & -\Delta^{-1}
\end{pmatrix}
$$
where $\Delta=P-N$ is a boundedly invertible selfadjoint positive definite operator on $X$. Moreover, 
\begin{equation}\label{diago_M}
T^{-1}MT = \begin{pmatrix}
A-BR^{-1}B^*N & 0 \\
0 & A-BR^{-1}B^*P
\end{pmatrix}
= \begin{pmatrix}
A_+ & 0 \\
0 & A_-
\end{pmatrix}
\end{equation}
by straightforward calculations. 
Moreover, subtracting the algebraic Riccati equations satisfied respectively by $P$ and $N$, we have
$$
\Delta A_- + A_+^* \Delta = 0 .
$$
This shows that the spectrum of $A_-$ is the negative of that of $A_+$.

\subsection{Consequence: proof of Proposition \ref{prop_optsol}}\label{sec_proof_exp}
Setting
\begin{equation}\label{def_vw}
\begin{pmatrix} v_T \\ w_T \end{pmatrix}
=
T^{-1} \begin{pmatrix} \delta y_T \\ \delta\lambda_T \end{pmatrix} ,
\end{equation}
we get from \eqref{extr_syst_delta} and \eqref{diago_M} that
\begin{equation*}
\begin{split}
\dot v_T &= A_+ v_T \\
\dot w_T &= A_- w_T \\
\end{split}
\end{equation*}
with $v_T(0)+w_T(0) = y_0-\bar y$ and $v_T(T)+w_T(T) = y_1-\bar y$.
Since
$v_T(t) 
= e^{-(T-t)A_+}v_T(T)$ and 
$w_T(t) = e^{tA_-}w_T(0) 
$,
taking $t=0$ yields
\begin{equation}\label{defUT}
\begin{pmatrix}
y_0-\bar y \\
y_1-\bar y 
\end{pmatrix}
=
\underbrace{
\begin{pmatrix}
\mathrm{id} & e^{-TA_+} \\
e^{TA_-} & \mathrm{id}
\end{pmatrix}
}_{J_T}
\begin{pmatrix}
w_T(0) \\
v_T(T)
\end{pmatrix} .
\end{equation}
By the exponential stability property, there exist $M>0$ and $\nu>0$ such that, using the operator norm, $\Vert e^{sA_-}\Vert_{L(X)}\leq Me^{-\nu s}$ and $\Vert e^{-sA_+}\Vert_{L(X)}\leq Me^{-\nu s}$ for every $s\geq 0$.
In particular, $\Vert e^{TA_-}\Vert_{L(X)}\leq Me^{-\nu T}$ and $\Vert e^{-TA_+}\Vert_{L(X)}\leq Me^{-\nu T}$. Therefore there exists $T_0>0$ such that, for every $T>T_0$, the operator $J_T$ (on $X$) defined in \eqref{defUT} is boundedly invertible, and we have
$$
J_T = \begin{pmatrix}
\mathrm{id} & 0 \\
0 & \mathrm{id}
\end{pmatrix}
+ \mathrm{O}(e^{-\nu T})
$$
and
\begin{equation*}
\begin{split}
J_T^{-1} 
&= 
\begin{pmatrix}
\mathrm{id} + e^{-TA_+} \left(\mathrm{id} + e^{TA_-} e^{-TA_+}\right)^{-1} e^{TA_-} & - e^{-TA_+} \left(\mathrm{id} + e^{TA_-} e^{-TA_+}\right)^{-1} \\
- \left(\mathrm{id} + e^{TA_-} e^{-TA_+}\right)^{-1} e^{TA_-} & \left(\mathrm{id} + e^{TA_-} e^{-TA_+}\right)^{-1}
\end{pmatrix} \\
&= \begin{pmatrix}
\mathrm{id} & -e^{-TA_+} \\
-e^{TA_-} & \mathrm{id}
\end{pmatrix}
+ \mathrm{O}(e^{-2\nu T})
\end{split}
\end{equation*}
where the remainder terms in $\mathrm{O}(\cdot)$ are in the sense of the operator norm as $T\rightarrow +\infty$.
Hence
\begin{equation*}
\begin{split}
w_T(0) &= y_0-\bar y - e^{-TA_+}(y_1-\bar y) + \mathrm{O}(e^{-2\nu T}) \\
v_T(T) &= - e^{TA_-}(y_0-\bar y) + y_1-\bar y + \mathrm{O}(e^{-2\nu T})
\end{split}
\end{equation*}
and thus
\begin{equation*}
\begin{split}
v_T(t) &= - e^{-(T-t)A_+} e^{TA_-}(y_0-\bar y) + e^{-(T-t)A_+} (y_1-\bar y) + e^{-(T-t)A_+} \mathrm{O}(e^{-2\nu T}) \\ 
w_T(t) &= e^{tA_-}(y_0-\bar y) - e^{tA_-} e^{-TA_+}(y_1-\bar y) + e^{tA_-}  \mathrm{O}(e^{-2\nu T})
\end{split}
\end{equation*}
as $T\rightarrow +\infty$.

Using \eqref{def_vw}, the formulas \eqref{exp_yT_1}, \eqref{exp_yT_2}, \eqref{exp_lambdaT_1}, \eqref{exp_lambdaT_2}, \eqref{exp_uT_1} and \eqref{exp_uT_2} of Proposition \ref{prop_optsol} follow.
The expansions can actually be obtained at any order, by expanding $J_T^{-1}$ to higher orders.

\subsection{Proof of Theorem \ref{main_thm}}\label{sec_proof_main_thm}
We infer from \eqref{exp_yT_1}, \eqref{zfzb} and \eqref{exp_uT_1} that
\begin{eqnarray}
V_T(y_0,y_1) &=& C_T(u_T) = \frac{1}{2}\int_0^T \left( \Vert y_T(t)-y_d\Vert_Q^2 + \Vert u_T(t)-u_d\Vert_R^2\right) dt \nonumber\\
&=& \frac{T}{2}\left( \Vert \bar y-y_d\Vert_Q^2 + \Vert \bar u-u_d\Vert_R^2\right) \label{expp1} \\
&& + \frac{1}{2}\int_0^T \left( \Vert z_f(t)\Vert_Q^2 + \Vert v_f(t)\Vert_R^2\right) dt \label{expp2}\\
&& + \frac{1}{2}\int_0^T \left( \Vert z_b(T-t)\Vert_Q^2 + \Vert v_b(T-t)\Vert_R^2\right) dt \label{expp3}\\
&& + \int_0^T \left( \langle Q(\bar y-y_d),z_f(t)\rangle_X + \langle R(\bar u-u_d),v_f(t)\rangle_U\right) \label{expp4}\\
&& + \int_0^T \left( \langle Q(\bar y-y_d),z_b(T-t)\rangle_X + \langle R(\bar u-u_d),v_b(T-t)\rangle_U\right) \label{expp4bis}\\
&& + \int_0^T \left( \langle Qz_f(t),z_b(T-t)\rangle_X + \langle Rv_f(t),v_b(T-t)\rangle_U\right) \label{expp5}\\
&& + \mathrm{O}(e^{-\nu T}) \label{expp6}
\end{eqnarray}
as $T\rightarrow+\infty$.
Above, the remainder term \eqref{expp6} is obtained by integration. The first term \eqref{expp1} is identified with $T\bar V$, where we recall that $\bar V$ is defined by \eqref{Vbar}. The second term \eqref{expp2} is equal to
$$
S_f(y_0) - \frac{1}{2}\int_T^{+\infty} \left( \Vert z_f(t)\Vert_Q^2 + \Vert v_f(t)\Vert_R^2\right) dt = S_f(y_0) + \mathrm{O}(e^{-\nu T})
$$
where we recall that $S_f(y_0)$ is the optimal value of the forward stabilization problem defined in Section \ref{sec_forward_stab}.
Similarly, the third term \eqref{expp3} is equal to $S_b(y_1) + \mathrm{O}(e^{-\nu T})$.

Let us compute the term \eqref{expp4}. Using \eqref{zfzb}, we have
$$
\int_0^T z_f(t)\, dt = A_-^{-1}(e^{TA_-}-\mathrm{id})(y_0-\bar y) = -A_-^{-1}(y_0-\bar y) + \mathrm{O}(e^{-\nu T}),
$$
and since $v_f(t)=-R^{-1}B^*Pz_f(t)$, we have
$$
\int_0^T v_f(t)\, dt = R^{-1}B^*PA_-^{-1}(y_0-\bar y) + \mathrm{O}(e^{-\nu T}).
$$
Then, using \eqref{KKT}, we infer that the term \eqref{expp4} is equal to
$$
\langle\bar\lambda,\underbrace{(A-BR^{-1}B^*P)}_{A_-}\int_0^Tz_f(t)\, dt\rangle_X
= \langle\bar\lambda,\bar y-y_0\rangle_X + \mathrm{O}(e^{-\nu T}).
$$
Similarly, the term \eqref{expp4bis} is equal to
$\langle\bar\lambda,y_1-\bar y\rangle_X + \mathrm{O}(e^{-\nu T})$.

Finally, using \eqref{zfzb}, the term \eqref{expp5} is a $\mathrm{O}(e^{-\nu T})$ as $T\rightarrow+\infty$.

We have thus obtained \eqref{asympt_exp}.

\subsection{Proof of Proposition \ref{prop_dissip}}
We can arrange differently the terms in the computations of the previous section. Gathering the terms \eqref{expp2} and \eqref{expp4} gives
\begin{equation*}
\begin{split}
& \frac{1}{2}\int_0^T \left( \Vert z_f(t)\Vert_Q^2 + \Vert v_f(t)\Vert_R^2 + 2 \langle z_f(t),Q(\bar y-y_d)\rangle_X + 2 \langle v_f(t),R(\bar u-u_d)\rangle_U \right) dt \\
=&\ \frac{1}{2}\int_0^T \Big( \Vert \underbrace{z_f(t)}_{y_f(t)-\bar y} + \bar y-y_d\Vert_Q^2 + \Vert \underbrace{v_f(t)}_{u_f(t)-\bar u}+\bar u-u_d\Vert_R^2 - \Vert \bar y-y_d\Vert_Q^2-\Vert\bar u-u_d\Vert_R^2 \Big) dt \\
=&\ \frac{1}{2}\int_0^T\left \Vert y_f(t)-y_d\Vert_Q^2 + \Vert u_f(t)-u_d\Vert_R^2 - \Vert \bar y-y_d\Vert_Q^2-\Vert\bar u-u_d\Vert_R^2 \right) dt \\
=&\ \frac{1}{2}\int_0^{+\infty}\left \Vert y_f(t)-y_d\Vert_Q^2 + \Vert u_f(t)-u_d\Vert_R^2 - \Vert \bar y-y_d\Vert_Q^2-\Vert\bar u-u_d\Vert_R^2 \right) dt + \mathrm{O}(e^{-\nu T}) \\
=&\ V_f(y_0) + \mathrm{O}(e^{-\nu T})
\end{split}
\end{equation*}
where $V_f(y_0)$ is defined by \eqref{def_Vf}. 
Similarly, the sum of the two terms \eqref{expp3} and \eqref{expp4bis} is equal to $V_b(y_1) + \mathrm{O}(e^{-\nu T})$. We have thus proved \eqref{asympt_exp2} and all contents of Proposition \ref{prop_dissip}.

\section{An additional result for free final states}
\subsection{Setting and main result}
In this section, we consider the optimal control problem \eqref{contsyst}-\eqref{cost} with fixed initial state $y(0)=y_0$ but with free final state $y(T)\in X$. 
The value function of such an optimal control problem is then defined by
$$
V_T(y_0) = \min\{ C_T(u)\ \mid\ y(0)=y_0\} .
$$
As before, there exists a unique solution $(y_T(\cdot),u_T(\cdot))$.
The objective is to provide an asymptotic expansion of $V_T(y_0)$ as $T\rightarrow+\infty$.

Compared with Section \ref{sec_main_result}, interestingly, we can relax the assumptions of exact controllability and of group generation and even consider unbounded control operators (occuring in boundary control problems). Also, we only assume that $Q$ is positive semidefinite. Hereafter, we assume that:
\begin{itemize}
\item The operator $A:D(A)\rightarrow X$ generates on $X$ an analytic $C_0$ semigroup $(e^{tA})_{t\geq 0}$.
\item The control operator $B\in L(X,D(A^*)')$ is admissible.
\item The pair $(A,B)$ is exponentially stabilizable and the pair $(A,Q^{1/2})$ is exponentially detectable.
\end{itemize}
The admissibility assumption implies that the control operator $B$ is 'not too unbounded' (see \cite{trelat_book, TucsnakWeiss} for details and examples). For instance, this framework includes the heat equation with Neumann boundary control. Refer to \cite{CurtainZwart} for concepts related to exponential stabilizability and detectability.   We  refer the reader to \cite{CurtainZwart} for concepts related to  exponential stabilizability and detectability.

\begin{theorem}\label{thm2}
Defining $V_f(y_0)$ as in \eqref{def_Vf}, there exists $\mu\in\R$ such that
\begin{equation}\label{expVT2}
V_T(y_0) = T\bar V + V_f(y_0) + \mu + \mathrm{O}(e^{-\nu T}) 
\end{equation}
as $T\rightarrow+\infty$. 
\end{theorem}

Additionally to the statement of Theorem \ref{thm2}, we have the following result.

\begin{proposition}\label{prop_optsol2}
We have
\begin{equation*}
\begin{split}
y_T(t) &= \bar y + e^{tA_-}(y_0-\bar y) - E e^{(T-t)A_-^*} w_T(T) + \mathrm{O}(e^{-\nu(t+T)}) \\
\lambda_T(t) &= \bar\lambda - P e^{tA_-}(y_0-\bar y) + (\mathrm{id}+PE) e^{(T-t)A_-^*} w_T(T) + \mathrm{O}(e^{-\nu(t+T)}) \\
u_T(t) &= \bar u - R^{-1}B^*Pe^{tA_-}(y_0-\bar y) + R^{-1}B^*(\mathrm{id}+PE) e^{(T-t)A_-^*} w_T(T) + \mathrm{O}(e^{-\nu(t+T)})
\end{split}
\end{equation*}
where $w_T(T)=P(y_T(T)-\bar y)-\bar\lambda$.
\end{proposition}

\begin{remark}\label{rem_mu}
Theorem \ref{main_thm} generalizes \cite{EKPZ} where an expansion similar to \eqref{asympt_exp} was derived in the finite-dimensional context.
The constant $\mu$ (denoted $\lambda$ in \cite{EKPZ}), which may be alternatively defined as the limit of $V_T(\bar y)-T\bar V$ as $T\rightarrow+\infty$, is related to ergodic considerations on the Hamilton-Jacobi equation -- here, coinciding with the algebraic Riccati equation (see, e.g., \cite{Arisawa, BarlesSouganidis, Ishii, QuincampoixRenault}) and is sometimes called the ergodic constant. 
We refer the reader to the numerous comments and citations done in \cite{EKPZ}.

Here, in addition to the extension to infinite dimension, we provide an expression for $\mu$:
\begin{multline}\label{mu}
\mu = \lim_{T\rightarrow+\infty} \frac{1}{2}\int_0^T \left( \Vert E e^{(T-t)A_-^*}w_T(T)\Vert_Q^2 + \Vert R^{-1}B^*(\mathrm{id}+PE)e^{(T-t)A_-^*}w_T(T)\Vert_R^2\right) dt \\
+ \Big\langle\bar\lambda, (-AE+BR^{-1}B^*(\mathrm{id}+PE))\int_0^Te^{(T-t)A_-^*}\, dt \, w_T(T)\Big\rangle_X
\end{multline}
where 
$E = -\int_0^{+\infty} e^{tA_-}BR^{-1}B^*e^{tA_-}\, dt$ is the unique solution of the Lyapunov equation $A_-E+EA_-^* - BR^{-1}B^* = 0$.
\end{remark}

\subsection{Proof of Theorem \ref{thm2}}\label{sec_proof_thm2}
The strategy is quite similar to the proof of Theorem \ref{main_thm} and actually, except some minor modifications (such as, in the PMP, we have $\lambda_T(T)=0$ by transversality, because $y_T(T)$ is let free), the proof is the same in finite dimension. But, in infinite dimension, the main difference is that now we have a semigroup, and not a group; as a consequence, $P>0$ still exists (but $P$ is not boundedly invertible in general) but the existence of $N$ is not ensured. Therefore, the strategy of diagonalization by blocks of $M$, developed in Section \ref{sec_diago}, needs to be adapted.
Actually, an appropriate strategy has been developed in \cite[Section 3.3, proof of Theorem 6]{TrelatZhangZuazua_SICON2018}. 
Let us sum up hereafter, rapidly, the main steps.

First of all, under the assumptions that have been done (admissibility of $B$, analyticity of the semigroup, exponential stabilizability and detectability), there exists a positive definite selfadjoint operator $P\in L(X)$ solution of the algebraic Riccati equation \eqref{ARE} (see \cite{LasieckaTriggiani, Flandoli}); note however that $P$ fails in general to be boundedly invertible. The operator $A_-=A-BR^{-1}B^*P$ generates on $X$ an exponentially stable $C_0$ semigroup $(e^{tA_-})_{t\geq 0}$. Hence, the forward stabilization problem considered in Section \ref{sec_forward_stab} is still well defined, with the same solution. Note that, to define properly $A_-$, it is required to consider the operator $B^*P$; since $B^*\in L(D(A^*),U)$, it could happen that this operator does not make sense. But it is part of the results contained in \cite{LasieckaTriggiani} that $B^*P\in L(X,U)$ (and this is a nontrivial issue). 

As a second step, it is proved in \cite{TrelatZhangZuazua_SICON2018} that there exists an operator $E\in L(X)$ solution of the Lyapunov equation
$$
(\underbrace{A-BR^{-1}B^*P}_{A_-})E+E(\underbrace{A-BR^{-1}B^*P}_{A_-})^* - BR^{-1}B^* = 0 .
$$
Actually, $E = -\int_0^{+\infty} e^{tA_-}BR^{-1}B^*e^{tA_-}\, dt$.
Now, defining
$$
L = \begin{pmatrix}
\mathrm{id}+EP & E \\
P & \mathrm{id}
\end{pmatrix} ,
$$
the operator $L\in L(X)$ is boundedly invertible and
$$
L^{-1} = \begin{pmatrix}
\mathrm{id} & -E \\
-P & \mathrm{id}+PE
\end{pmatrix} ,
$$
and we have
$$
LML^{-1} = \begin{pmatrix} A_- & 0 \\ 0 & -A_-^* \end{pmatrix} .
$$
Then, setting
$$
\begin{pmatrix}
v_T \\ w_T
\end{pmatrix}
=
L
\begin{pmatrix}
\delta y_T \\ \delta\lambda_T
\end{pmatrix} ,
$$
we have
\begin{equation*}
\begin{split}
\dot v_T &= A_- v_T \\
\dot w_T &= -A_-^* w_T
\end{split}
\end{equation*}
with $v_T(0)-Ew_T(0)=y_0-\bar y$. 
We have $v_T(t)=e^{tA_-}v_T(0)$ and $w_T(t)=e^{(T-t)A_-^*}w_T(T)$.
By \cite[Lemma 2]{TrelatZhangZuazua_SICON2018}, there exists $C>0$ (not depending on $T$) such that $\Vert\delta y_T(T)\Vert_X+\Vert\delta\lambda_T(0)\Vert_X\leq C(\Vert\delta y(0)\Vert_X+\Vert\delta\lambda_T(T)\Vert_X)$. Since $\delta y_T(0)=y_0-\bar y$ and $\delta\lambda_T(T)=-\bar\lambda$, it follows that $\delta\lambda_T(0)$ and $\delta y_T(T)$ are bounded in $X$ uniformly with respect to $T$, and thus $v_T(0)$ and $w_T(T)$ are bounded in $X$ uniformly with respect to $T$. Therefore $v_T(T)=\mathrm{O}(e^{-\nu T})$ and $w_T(0)=\mathrm{O}(e^{-\nu T})$.
As a consequence,
$$
v_T(0) = y_0-\bar y + \mathrm{O}(e^{-\nu T}), \qquad
(\mathrm{id}+PE)w_T(T) = -\bar\lambda + \mathrm{O}(e^{-\nu T}) .
$$
We do not know whether $\mathrm{id}+PE$ is boundedly invertible or not, but we insist that we know that $w_T(T)$ is bounded in $X$ uniformly with respect to $T$. Finally, we obtain
$$
v_T(t) = e^{tA_-}(y_0-\bar y) + e^{tA_-} \mathrm{O}(e^{-\nu T}),
\qquad
w_T(t) = e^{(T-t)A_-^*} w_T(T) .
$$
Note that, since $w_T=P(y_T-\bar y)+\lambda_T-\bar\lambda$ and since $\lambda_T(T)=0$ (transversality condition in the PMP, since $y_T(T)$ is let free), we have $w_T(T)=P(y_T(T)-\bar y)-\bar\lambda$.
We then infer the results of Proposition \ref{prop_optsol2}.

To establish the statements of Theorem \ref{thm2} and Remark \ref{rem_mu}, we proceed as in Section \ref{sec_proof_main_thm}, by gathering adequately the various terms. We do not give any details.
The fact that $\mu$ is well defined by \eqref{mu} is because we first define $\mu$ as the limit of $V_T(y_0)-T\bar V-V_f(y_0)$ as $T\rightarrow+\infty$ (as it is the usual definition for the ergodic constant) and we identify this limit in the exponential scale provided by the asymptotic expansion.

\section{Conclusion}
For large-time LQ optimal control problems with fixed terminal points, we presented an asymptotic expansion of the value function. The first term is obtained using the turnpike property, and the second term is the sum of the optimal values of two stabilization problems corresponding to the respective terminal points towards the turnpike. We also established a version of this result when the final point is left free.

Furthermore, in both cases, we derived explicit expansions for the optimal trajectories. Theorems \ref{main_thm} and \ref{thm2} are applicable in infinite-dimensional spaces but under different assumptions. In Theorem \ref{main_thm}, we assumed exact controllability, group generation, and a bounded control operator. On the other hand, Theorem \ref{thm2} is based on assumptions of an analytic semigroup, a possibly unbounded but admissible control operator, and exponential stabilizability and detectability.

These assumptions enable us to use the comprehensive algebraic Riccati theory, akin to the one in finite dimension, which has been instrumental in proving the theorems.

\paragraph{Relaxing analyticity}
If one wishes to relax the analyticity assumption, the algebraic Riccati theory becomes considerably more complicated. In particular, there exist several possible algebraic Riccati equations with different meanings and interpretations (see \cite{FlandoliLasieckaTriggiani, WeissRebarber, WeissZwart}). Consequently, it is not clear whether Theorem \ref{thm2} and its proof, developed in Section \ref{sec_proof_thm2}, can be adapted to this more general context. We leave this issue open for further exploration.

To be more precise, we refer to the explanation given in \cite[Section 3.3, page 108]{TrelatWangXu}, that we reproduce partly here. The operator $A_-=A-BR^{-1}B*P$, of domain $D(A_-)$ (which may differ from $D(A)$) generates an exponentially stable $C_0$ semigroup $(e^{tA_-})_{t\geq 0}$. Here, $P\in L(X)$ is a positive definite selfadjoint operator that maps $D(A_-)$ to $D(A^*)$ and $D(A)$ to $D(A_-^*)$, and we have $S_f(y_0) = \frac{1}{2}\langle P(y_0-\bar y),y_0-\bar y\rangle_X$. However, $P$ satisfies an algebraic Riccati on $D(A_-)$ and possibly another one on $D(A)$. It is not known whether the control operator $B$ is admissible or not for the semigroup $(e^{tA_-})_{t\geq 0}$. The latter fact is a serious obstacle to defining $E = -\int_0^{+\infty} e^{tA_-}BR^{-1}B^*e^{tA_-}\, dt$ solution of the Lyapunov equation $A_-E+EA_--BR^{-1}B^*=0$, and thus to defining the operator $L$ instrumentally used in Section \ref{sec_proof_thm2} to diagonalize $M$.

\paragraph{Nonlinear systems}
Another final remark is the following. In the present paper we have treated LQ optimal control problems. It is natural to ask whether our results can be established in a nonlinear context, for the optimal control problem
\begin{equation*}
\begin{split}
& \dot y(t) = f(y(t),u(t)),\qquad y(0)=y_0,\quad y(T)=y_1, \\
& V_T(y_0,y_1) = \min C_T(u)\qquad\textrm{where}\qquad C_T(u) = \int_0^T f^0(y(t),u(t))\, dt .
\end{split}
\end{equation*}
This study is done in \cite{ATZ} and since it is interesting to compare the obtained results with those of the present paper we provide hereafter a quick overview of the main result of \cite{ATZ}.

Let $(\bar y,\bar u)$ (the turnpike) be a solution of the static optimal control problem
$$
\bar V = \min_{f(y,u)=0}f^0(y,u) .
$$
We define the forward stabilization problem as
\begin{equation*}
\begin{split}
& \dot y(t)=f(y(t),u(t)),\qquad y(0)=y_0, \\
& V_f(y_0) = \min\int_0^{+\infty}\left( f^0(y(t),u(t))-f^0(\bar y,\bar u)\right) dt
\end{split}
\end{equation*}
and the backward stabilization problem as
\begin{equation*}
\begin{split}
& \dot y(t)=-f(y(t),u(t)),\qquad y(0)=y_1, \\
& V_b(y_1) = \min\int_0^{+\infty}\left( f^0(y(t),u(t))-f^0(\bar y,\bar u)\right) dt .
\end{split}
\end{equation*}
Noting that $w(y,u)=f^0(y,u)-f^0(\bar y,\bar u)$ is the usual supply rate function that is used to characterize the dissipativity property of an optimal control problem (see \cite{Willems} for the notion of dissipativity and see \cite{Faulwasser, GruneGuglielmi, GruneMuller, TrelatZhang} for the various relationships between dissipativity and turnpike). 
It is proved in \cite{ATZ} that, in finite dimension, under dissipativity and other appropriate assumptions,
$$
V_T(y_0,y_1) = T\bar V + V_f(y_0) + V_b(y_1) + \mathrm{o}(1)
$$
as $T\rightarrow+\infty$. 

It is interesting to compare this result to the results obtained in the present paper -- and more precisely, with the contents of Proposition \ref{prop_dissip}, which are more prepared for the comparison with respect to the dissipativity property, as already alluded. Apart from the dimension, the results in the LQ case stated in Propositions \ref{prop_dissip} and \ref{prop_optsol} are more precise, in that, not only, the remainder term is not a $\mathrm{o}(1)$ and a $\mathrm{O}(e^{-\nu T})$, but also, in the LQ case we have obtained an expansion for the optimal solution, which we do not have in the nonlinear case.

\end{document}